\documentclass{article}

\usepackage[proportional,scaled=1.064]{erewhon}
\usepackage[T1]{fontenc}
\usepackage{mathrsfs}

\usepackage{geometry}
\selectfont

\usepackage[utf8]{inputenc}
\usepackage[dvipsinames]{xcolor}
\usepackage[absolute,overlay]{textpos}
\usepackage{graphicx} 
\usepackage{array}
\usepackage{hyperref}

\usepackage{natbib}
\usepackage[utf8]{inputenc} 
\usepackage[T1]{fontenc}    
\usepackage{url}            
\usepackage{booktabs}       
\usepackage{amsfonts}       
\usepackage{nicefrac}       
\usepackage{dsfont}            
\usepackage{import}
\usepackage[toc,page]{appendix}
\usepackage{minitoc}

\usepackage{amsmath}
\usepackage{amssymb,amsthm}

\usepackage{mathrsfs}
\usepackage{nicefrac}

\definecolor{darkmidnightblue}{HTML}{003366}    
\definecolor{midnightblue}{HTML}{0059b3}
\definecolor{chromered}{HTML}{f14233}

\definecolor{darkpowderblue}{rgb}{0.0, 0.2, 0.6}
\definecolor{dukeblue}{rgb}{0.0, 0.0, 0.61}

\hypersetup{
    colorlinks = true,
    citecolor= darkmidnightblue,
    urlcolor=darkmidnightblue,
    breaklinks=true,
    linkcolor = darkmidnightblue,
    linkbordercolor = {white},
}

\usepackage{cleveref}
\usepackage{autonum}
\usepackage{algorithm}
\usepackage{algorithmic}

\usepackage[absolute,overlay]{textpos}

\newcommand\bs{\boldsymbol}

\newcommand\bfE{\mathbf E}

\newcommand\bY{\boldsymbol Y}

\newcommand\bu{\boldsymbol u}
\newcommand\bv{\boldsymbol v}

\newcommand\bxi{\boldsymbol\xi}

\newcommand\bzeta{\bs \zeta}

\newcommand\btheta{\boldsymbol\theta}

\def\hat{\widehat}

\def\bfE{\mathbf E}

\newcommand{\textupsf}[1]{\textup{\textsf{#1}}}
\def\CLPEN1{\hyperlink{PEN}{\textupsf{{$(C,\lambda)$-PEN}}}}

\def\PLD1{\hyperlink{PLD}{\textupsf{{PLD}}}}
\def\PGF1{\hyperlink{PGF}{\textupsf{{PGF}}}}
\def\TPLD1{\hyperlink{TPLD}{\textupsf{{TPLD}}}}

\newtheorem{theorem}{Theorem}

\newtheorem{remark}{Remark}

\crefname{cor}{Corollary}{Corollaries}
\crefname{condition}{Condition}{C}

\def\bu{\boldsymbol u}

\def\btheta{\boldsymbol\theta}
\def\balign#1\ealign{\begin{align}#1\end{align}}
\def\baligns#1\ealigns{\begin{align*}#1\end{align*}}
\def\balignat#1\ealign{\begin{alignat}#1\end{alignat}}
\def\balignats#1\ealigns{\begin{alignat*}#1\end{alignat*}}
\def\bitemize#1\eitemize{\begin{itemize}#1\end{itemize}}
\def\benumerate#1\eenumerate{\begin{enumerate}#1\end{enumerate}}
\newenvironment{talign}
 {\align}
 {\endalign}
\newenvironment{talign*}
 {\csname align*\endcsname}
 {\endalign}

\theoremstyle{definition}
\newtheorem{corollary}{Corollary}

\title{Simple proof of the risk bound for denoising by 
exponential weights for asymmetric noise distributions}

\author{Arnak S. Dalalyan\\[5pt] CREST, ENSAE, Institut Polytechnique de Paris}

\date{\vspace{-5ex}}

\begin{document}

\maketitle

\begin{abstract}
    In this note, we consider the problem of aggregation of
    estimators in order to denoise a signal. The main 
    contribution is a short proof of the fact that the
    exponentially weighted aggregate satisfies a sharp 
    oracle inequality. While this result was already known 
    for a wide class of symmetric noise distributions, the 
    extension to asymmetric distributions presented in 
    this note is new.
\end{abstract}

\section{Introduction}

Let us consider the problem of denoising an $n$ dimensional
noisy signal $\bY$ using a family of candidates $\btheta_1,
\ldots,\btheta_m$. More precisely, we assume that
\begin{align}
    \bY = \btheta^* + \bxi
\end{align}
where $\btheta^*\in\mathbb R^n$ is the $n$ dimensional 
true signal and $\bxi$ is random noise. Only the noisy 
vector $\bY$ is observed and the goal is to construct 
an estimator $\hat\btheta$ such that the expected error 
$\bfE[\|\hat\btheta -\btheta^*\|^2]$ is as small as 
possible, where $\|\bv\|$ stands for the Euclidean 
norm of $\bv\in\mathbb R^n$. We consider the framework 
in which to achieve the aforementioned goal we are 
given a set of vectors $\{\btheta_1,\ldots,\btheta_m\}$. 
An estimator $\hat\btheta$ is considered a good estimator, 
if the regret
\begin{align}\label{eq:regret}
    \bfE[\|\hat\btheta -\btheta^*\|^2] - \min_{j=1,
    \ldots,m}\|\btheta_j -\btheta^*\|^2
\end{align}
is as small as possible. This problem has been coined
model-selection aggregation in \citep{Tsyb03}, where it
also proved that the optimal rate of the difference  
\eqref{eq:regret} is $\log m$. The problem of aggregation
has been extensively studied in the literature, see
for instance \citep{BTW2007,Yang00,Yang1,Yang2,JRT,bellec2018,Rigollet11,TsybICM,Alquier3, Lecue,Golubev1}. In this note, we consider the
exponentially weighted aggregate (EWA) defined as follows.
Let $\pi_0(1),\ldots,\pi_0(m)$ be some nonnegative weights 
summing to one. Each $\pi_0(j)$ represent our prior 
confidence in the approximation of $\btheta^*$ be 
$\btheta_j$. Based on these prior weights and the observed 
vector $\bY$, we define
\begin{align}
    \hat\btheta = \sum_{j=1}^m \btheta_j \hat\pi(j),
    \qquad\text{with}\qquad \hat\pi(j) = \frac{\exp\{ -
    \|\bY-\btheta_j\|^2/\beta\}\pi_0(j)}{\sum_{\ell=1}^m 
    \exp\{ -\|\bY-\btheta_\ell\|^2/\beta\}\pi_0(\ell)}.
\end{align}
In this expression, $\beta>0$ is a tuning parameter of the
method. As established in the aforementioned references, in
different settings one can prove that EWA satisfies the
inequality
\begin{align}\label{EWA1}
    \bfE[\|\hat\btheta -\btheta^*\|^2] \leqslant \min_{j=1,
    \ldots,m}\Big(\|\btheta_j -\btheta^*\|^2 + \beta
    \log(1/\pi_0(j))\Big).
\end{align}
In particular, if $\pi_0$ is the uniform distribution over
$\{1,\ldots,m\}$, one obtains the rate-optimal remainder
term $\beta\log m$ for the difference in \eqref{eq:regret}.

As pointed out in some papers \citep{DT07,DT08,Dal_IHP}, 
it is helpful to extend the above-described framework to
the case of aggregating a family of estimators which is
potentially infinite. This is equivalent to considering
a subset $S_0\subset \mathbb R^n$ and aiming at finding an
``optimal'' way of combining all its elements in order
to estimate $\btheta^*$. These types of considerations have
led to the following extension of the estimator \eqref{EWA1}:
\begin{align}\label{EWA2}
    \hat\btheta = \int_{\mathbb R^n} \btheta\,\hat\pi(d\btheta),
    \qquad \text{with}\qquad \frac{d\hat\pi}{d\pi_0} 
    (\btheta) = \frac{\exp\{-\|\bY-\btheta\|^2/\beta\}}{
    \int_{\mathbb R^n} \exp\{-\|\bY-\bu\|^2/\beta\}\pi_0(d\bu)
    } .
\end{align}
Notice that this estimator is the Bayesian posterior mean
in the case where $\bxi$ is drawn from the Gaussian 
distribution with zero mean and covariance matrix 
$(\beta/2)\mathbf I_n$. The goal of this note is to 
provide an alternative and simple proof of the fact that
EWA $\hat\btheta$ satisfies \eqref{EWA1} and its extension
to aggregating an infinite set, provided that the 
distribution of the noise $\xi$ satisfies some suitable 
conditions. We also slightly extend the existing results by
including noise distributions that are not symmetric with
respect to the origin. This is particularly suitable for
estimating the parameters of Bernoulli or binomial 
distributions. 

\paragraph{Notation} We use boldface letters for 
vectors, which are always seen as one-column matrices.  
For any vector $\bv$, $\|\bv\|$ and $\|\bv\|_\infty$ 
are respectively the Euclidean norm and the sup-norm. 
By convention, throughout this work, $0\cdot\infty = 0$. 
For a probability distribution $\pi$ on $\mathbb R^n$, we 
denote by $\text{Var}_\pi(\btheta)$ the variance with 
respect to $\pi$ defined by $\int_{\mathbb R^n} 
\|\btheta\|^2\,\pi(d\btheta) - \|\int_{\mathbb R^n} 
\btheta\,\pi(d\btheta)\|^2$. For two probability 
distributions $\mu$ and $\nu$ defined on the same 
probability space and such that $\mu$ is absolutely 
continuous with respect to $\nu$, the Kullback-Leibler
divergence is defined by $D_{\rm KL} (\mu||\nu) =  
\int \frac{d\mu}{d\nu}(x)\log \frac{d\mu}{d\nu}(x)\,
\nu(dx)$.

\section{Main result}
This section is devoted to stating and briefly discussing
the main result, the proof being postponed to \Cref{sec:proof} 
below. Prior to stating the result, we recall the Bernstein condition. For some $v>0$ and $b\geqslant 0$, we 
say that a random variable $\eta$ satisfies the $(v,b) 
$-Bernstein condition, if 
\begin{align}
    \mathbf E[e^{t\eta}] \leqslant \exp\Big\{\frac{v^2t^2}{2
    (1-b|t|)}\Big\},\qquad \forall t\in(-1/b,1/b).
\end{align}
This condition is clearly on the distribution of the random variable.
One can check that if $\eta$ satisfies the $(v,b)$-Bernstein
condition, then it is sub-exponential with zero mean, and the 
variance of $\eta$ is at least equal to $v$. Many common 
distributions satisfy this assumption. For instance, any
sub-Gaussian distribution with variance proxy $\tau$ satisfies
the $(\tau,0)$-Bernstein condition. Any random variable supported
by $[-A,A]$ satisfies the Bernstein condition with $(v,b) = 
(A^2,0)$ but also with $(v,b) = (\text{Var}(\eta) , A/3)$ 
\citep{vershynin_2018}. We will see that the latter is more
useful for our purposes than the former. 

Similarly, if $\mathcal F$ is a sigma-algebra and $v$ and $b$ 
are two $\mathcal F$- measurable random variables, we say that 
$\eta$ is $(v,b)$-Bernstein conditionally to $\mathcal F$, if
almost surely, the inequality $\bfE[e^{t\eta}|\mathcal F] 
\leqslant \exp\{v^2t^2/(1-b|t|)\}$ is satisfied for every 
$t\in\mathbb R$ such that $|t|b<1$.

\begin{theorem}\label{thm}
    Let $\pi_0$ be a probability distribution supported
    by $S_0\subset\mathbb R^n$ with a diameter measured in 
    sup-norm bounded by $\mathcal D_0$. Assume that the 
    distribution of $\bxi$ satisfies the following assumption: 
    for some sigma algebra $\mathcal F$ and for some 
    $b:[0,1]\to [0,\infty)$ and continuously differentiable 
    function $v:[0,1]\to [0,\infty)$ vanishing at the 
    origin,  for every $\alpha \in(0,1]$, there exists 
    an  $n$-dimensional random vector $\bzeta$ such that
    \begin{align}\label{cond:zeta}
        \bfE[\bzeta|\mathcal F] = 0,\qquad \bxi + \bzeta 
        \stackrel{\mathscr D}{=} (1+\alpha) \bxi . 
    \end{align}
    and, conditionally to $\mathcal F$, the entries $\zeta_i$ 
    are independent and satisfy the $(v(\alpha),b(\alpha))$-
    Bernstein condition. Then, for every $\beta \geqslant 
    2b(0) \mathcal D_0$, we have
    \begin{align}
        \bfE[\|\hat\btheta-\btheta^*\|^2]
        \leqslant \inf_{\pi}\bigg\{\int_{\mathbb R^n} 
        \|\btheta-\btheta^*\|^2 \,\pi(d\btheta) + \beta 
        D_{\rm KL}(\pi||\pi_0) \bigg\} +\bigg(\frac{2
         v'(0)}{ \beta -2 b(0) \mathcal D_0} 
        - 1\bigg) \bfE[{\rm Var}_{\hat\pi}(\btheta)],
        \label{eq:main}
    \end{align}
    where the $\inf$ is over all the probability distributions. 
    As a consequence, for $\beta \geqslant  2v'(0) 
    + 2b(0) \mathcal D_0$, we get
    \begin{align}
        \bfE[\|\hat\btheta-\btheta^*\|^2]
        \leqslant \inf_{\pi}\bigg\{\int_{\mathbb R^n} 
        \|\btheta-\btheta^*\|^2 \,\pi(d\btheta) + \beta 
        D_{\rm KL}(\pi||\pi_0) \bigg\}.
        \label{eq:main1}
    \end{align}
\end{theorem}

Let us briefly comment on this result. First, the link 
between \eqref{eq:main1} and \eqref{EWA1} might not be
easy to see.  It is obtained by considering a prior
distribution $\pi_0$ supported by the finite set $\{\btheta_1, 
\ldots,\btheta_m\}$ and by upper bounding the infimum in
\eqref{eq:main1} by the minimum over all the Dirac measures
$\delta_{\btheta_j}$. One easily checks that $D_{KL}(
\delta_{\btheta_j}||\pi_0) = \log(1/\pi_0(j))$, which allows
to infer \eqref{EWA1} from \eqref{eq:main1}. 

Second, one may wonder where the form of the upper bound in
\eqref{eq:main1} comes from. The presence of the KL-divergence
in this bound may seem surprising. The reason is that 
there is a deep connection between the KL-divergence and the
exponential weights. Indeed, according to the Varadhan-Donsker variational formula, the ``posterior'' distribution $\hat\pi$ 
defined in \eqref{EWA2} is solution to following problem:
\begin{align}\label{eq:DV}
    \hat\pi \in\mathop{\rm argmin}_\pi\limits \bigg\{\int_{
    \mathbb R^n} \|\btheta-\bY\|^2 \,\pi(d\btheta) + \beta 
    D_{\rm KL}(\pi||\pi_0) \bigg\},
\end{align}
where the $\min$ is over all the probability distributions. 
This result will be the starting point of the proof.

Finally, one can wonder how restrictive the assumptions of this
theorem are. We will show below that they are satisfied for a 
broad class of noise distributions.

    \section{Instantiation to some well-known noise 
    distributions}

    The main theorem stated in the previous section requires
    a general and a rather abstract condition to be satisfied
    by the noise distribution. This section shows that many
    distributions encountered in applications satisfy this
    assumption with some parameters $v'(0)$ and $b(0)$ 
    which are easy to determine.

    \subsection{Centered Bernoulli noise}
    
    Assume that each $\xi_i$
    is a centered Bernoulli random variable: it takes the 
    value $1-\rho_i$ with probability $\rho_i$ and the value 
    $-\rho_i$ with probability $1-\rho_i$. Here, $\rho_i\in 
    (0,1)$. Then, one can set 
    \begin{align}
        \mathbf P\big(\zeta_i = \alpha\xi_i\,|\, \xi_i \big)  
        = \frac{1+\alpha-\alpha|\xi_i|}{\alpha+1}
        ,\quad \mathbf P\big(\zeta_i = -\mathop{\rm sgn}(\xi_i
        )(1+\alpha - \alpha|\xi_i|)\,|\, \xi_i\big)  = \frac{
        \alpha |\xi_i|}{\alpha+1}.
    \end{align}
    We see that conditionally to $\xi_i$, the random variable
    $\zeta_i$ is zero mean and takes its values in an interval 
    of length $\alpha(1-\rho_i) + \alpha\rho_i +1 = \alpha\rho_i 
    + 1 + \alpha - \alpha\rho_i = 1+\alpha$. This implies that 
    $\zeta_i$ satisfies the $((1+\alpha)^2/4,0)$-Bernstein
    condition, conditionally to $\xi_i$. In other terms, 
    $\zeta_i$ is sub-Gaussian with variance proxy $(1 + 
    \alpha)^2/4$. However, this does not help in applying 
    \Cref{thm}, since the function $v(\alpha) = (1+\alpha)^2/4$
    does not vanish at the origin. On the positive side, since the 
    conditional variance of $\zeta_i$ given $\xi_i$ is smaller than $\alpha(1+\alpha)$ and the support is included in $[-(1+\alpha
    ), (1+\alpha)]$, the conditional distribution of $\zeta_i$ 
    given $\xi_i$ satisfies the Bernstein condition with 
    $v(\alpha) = \alpha(1+\alpha)$ and $b(\alpha) = (1+\alpha)/3$,
    see  \citep[Exercise 2.8.5]{vershynin_2018}. This 
    yields the following result.
    \begin{corollary}
        Let $\pi_0$ be a probability distribution supported by
        $S_0\subset \mathbb R^n$ such that $\mathcal D_0 = 
        \sup_{\btheta,\btheta'\in S_0}\|\btheta-\btheta'
        \|_\infty <\infty$. Assume that $\bxi$ has 
        independent entries  $\xi_i$ satisfying $\mathbf 
        P(\xi_i = 1-\rho_i) = 1 - \mathbf P(\xi_i = -\rho_i) 
        = \rho_i$ for some $\rho_i\in (0,1)$. Then, for every 
        $\beta \geqslant (2/3)\mathcal D_0$, we have    
        \begin{align}
            \bfE[\|\hat\btheta-\btheta^*\|^2]
            \leqslant \inf_{\pi}\bigg\{\int_{\mathbb R^n} 
            \|\btheta-\btheta^*\|^2 \,\pi(d\btheta) + \beta 
            D_{\rm KL}(\pi||\pi_0) \bigg\} +\Big(\frac{6}{3 
            \beta - 2\mathcal D_0} - 1 \Big) \bfE[{\rm Var}_{
            \hat\pi} (\btheta)].
            \label{eq:bernoulli}
        \end{align}
        In particular, if $\beta\geqslant 2 + (2/3)\mathcal D_0$, 
        the last term in \eqref{eq:bernoulli} is nonpositive 
        and, therefore, can be neglected.
    \end{corollary}

    This corollary can be used in cases where the observations
    $Y_i$ are independent Bernoulli random variables with 
    mean $\theta_i^*$. In such a situation, it is natural to
    choose a prior distribution $\pi_0$ that is concentrated
    on the unit hypercube $[0,1]^n$, the diameter of which in
    sup-norm is equal to $1$. The corollary implies that in
    such a situation the inequality stated in \eqref{eq:main1}
    is true provided that $\beta \geqslant 8/3$. We refer the
    reader to \citep{Etienne1} for an application of this result 
    to graphon estimation. 

    \subsection{Gaussian noise}

    In the case of the Gaussian noise $\bxi$ with independent 
    entries having $0$ mean and variance equal to $\sigma_i^2$,
    one can check that the conditions of \Cref{thm} are 
    satisfied with the random vector $\bzeta$ which is 
    independent of $\bxi$ and has independent entries 
    drawn from the Gaussian distribution $\mathcal N(0, 
    (2\alpha+\alpha^2)\sigma_i^2)$. This means that in the 
    Bernstein condition one can choose $\mathcal F = \sigma 
    (\bxi)$, $b = 0$ and $v(\alpha) = (2\alpha + \alpha^2)
    \max_{1\leqslant i\leqslant n}\sigma_i^2$, which leads 
    to the following result. 

    \begin{corollary}
        Let $\pi_0$ be a probability distribution on $\mathbb 
        R^n$. Assume that $\bxi$ has independent entries 
        $\xi_i\sim \mathcal N(0,\sigma_i^2)$, $i=1,\ldots,n$. 
        Then, for every $\beta >0$, we have    
        \begin{align}
            \bfE[\|\hat\btheta-\btheta^*\|^2]
            \leqslant \inf_{\pi}\bigg\{\int_{\mathbb R^n} 
            \|\btheta-\btheta^*\|^2 \,\pi(d\btheta) + \beta 
            D_{\rm KL}(\pi||\pi_0) \bigg\} +\big(
            4\sigma^2\beta^{-1} - 1\big) \bfE[{\rm Var}_{
            \hat\pi}(\btheta)],
            \label{eq:Gaussian}
        \end{align}
        where $\sigma = \max_{1\leqslant i\leqslant n} \sigma_i$.
        In particular, if $\beta\geqslant 4\sigma^2$, the last 
        term in \eqref{eq:Gaussian} is nonpositive and, therefore,
        can be neglected.
    \end{corollary}

    Some preliminary versions of this result can be traced back 
    to \citep{George86,George86b}. In the form \eqref{eq:main1}, 
    and with an
    extension to aggregation of projection estimators, the result
    appeared in \citep{LeungBarron}. Further generalisations to 
    various families of linear estimators have been explored in 
    \citep{DalSal}. The proof of the oracle inequality in all 
    these papers is very specific to the Gaussian distribution
    since it is based on Stein's lemma (integration by parts for
    the Gaussian measure). The alternative proof presented in this
    work relies on techniques developed in 
    \citep{DT07,colt_DalalyanT09, Dal_IHP}. 

    \subsection{Bounded noise}

    For every $a,b>0$, let $\mathcal B(a,b)$ be the distribution
    of a random variable that takes the values $a$ and $-b$ with
    probabilities $b/(a+b)$ and $a/(a+b)$, respectively. 
    If the distribution of $\xi_i$ can be written as a mixture
    of the distributions $\mathcal B(a,b)$ with a mixing distribution
    with bounded support, then our main theorem can be applied. 
    More precisely, assume that the distribution of $\xi_i$ is 
    given by 
    \begin{align}
        p_{\xi_i}(dx) = \int_{0}^A\int_0^B \frac{b\delta_a(dx) + 
        a\delta_{-b}(dx)}{a+b}\, \nu_i(da,db),
    \end{align}
    where $\nu_i$ is a probability distribution on $[0,A]
    \times [0,B]$. 
    This means that $\xi_i = \eta_i^{\alpha_i,\beta_i}$ 
    with random variables $(\alpha_i,\beta_i)$ drawn from 
    $\nu_i$ and $\eta_i^{a,b}$ drawn from the binary distribution 
    $\frac{b\delta_a(dx) + a\delta_{-b}(dx)}{a+b}$. Akin to the
    first subsection of this section, one can choose $\zeta_i^{a,b}$
    so that $(1+\alpha)\eta_i^{a,b}$ has the same distribution
    as $\eta_i^{a,b} + \zeta_i^{a,b}$, for every pair $(a,b)$. 
    Then, clearly, $(1+\alpha)\xi_i$ has the same distribution as
    $\xi_i+\zeta_i^{\alpha,\beta}$. Let $\mathcal F$ be the 
    sigma algebra generated by the random variables $\alpha, 
    \beta,\{\xi_i^{a,b}:(a,b)\in[0,A]\times[0,B],i\in[n]\}$. 
    Conditionally to $\mathcal F$, $\zeta_i^{a,b}$ is a binary 
    random variable with zero mean and takes its values in the interval $[-B,A]$, 
    it satisfies the Bernstein condition with $b(\alpha) = 
    (A+B)(1+\alpha)/3$ and $v(\alpha) = (A+B)^2\alpha (1+\alpha)$. Therefore, we get the following result. 

    \begin{corollary}
        Let $\pi_0$ be a probability distribution supported by
        $S_0\subset \mathbb R^n$ such that $\mathcal D_0 = 
        \sup_{\btheta,\btheta'\in S_0}\|\btheta-\btheta'
        \|_\infty <\infty$. Assume that $\bxi$ has 
        independent entries $\xi_i$, $i=1,\ldots,n$, taking 
        values in an interval $I_i$ of length at most $L$. 
        Then, for every $\beta \geqslant (2/3)L\mathcal D_0$, 
        we have    
        \begin{align}
            \bfE[\|\hat\btheta-\btheta^*\|^2]
            \leqslant \inf_{\pi}\bigg\{\int_{\mathbb R^n} 
            \|\btheta-\btheta^*\|^2 \,\pi(d\btheta) + \beta 
            D_{\rm KL}(\pi||\pi_0) \bigg\} +\Big(\frac{
            6L^2}{3\beta - 2L\mathcal D_0} - 1\Big) 
            \bfE[{\rm Var}_{\hat\pi}(\btheta)].
            \label{eq:bounded}
        \end{align}
        In particular, if $\beta\geqslant 2L^2 + (2/3)L
        \mathcal D_0$, the last term in \eqref{eq:bounded} 
        is nonpositive and, therefore, can be neglected.
    \end{corollary}

    This result is well suited for the setting
    where the components $Y_i$ of the observation $\bY$ are 
    bounded. For instance, if we know that $\mathbf P(Y_i\in 
    [0,L]) = 1$ for every $i\in\{1,\ldots,n\}$, then it is also
    natural to choose a prior distribution satisfying $\mathcal 
    D_0 = L$. Inequality \eqref{eq:main1} is then satisfied for
    every $\beta\geqslant (8/3)L^2$. Note that, to the best of
    our knowledge, this is the first time that such a precise
    bound is obtained for asymmetric noise distributions. The
    similar result established in \citep[Theorem 2]{Dal_IHP} 
    deals with symmetric distributions only. 

    \subsection{Centered binomial noise}

    Consider the case where $\xi_i$'s are independent and 
    drawn from a centered and scaled binomial distribution
    $a\mathcal B(k,\rho_i) - a k \rho_i$, where $a>0$ is the
    scaling factor. This distribution is a particular case
    of distributions supported by a finite interval considered
    in the previous subsection. One can therefore apply the
    last corollary with $L = ak$. However, this leads to
    a bound which is too crude. Indeed, one can use the fact
    that $\xi_i $ is equal in distribution to $a(\eta_1+\ldots
    + \eta_k)$ where $\eta_j$'s are iid centered Bernoulli
    variables. Defining $\bar\zeta_1,\ldots,\bar\zeta_k$ 
    as independent random variables satisfying
    \begin{align}
        \mathbf P\big(\bar\zeta_j = \alpha\eta_j\,|\, \eta_j \big)  
        = \frac{1+\alpha-\alpha|\eta_j|}{\alpha+1}
        ,\quad \mathbf P\big(\bar\zeta_j = -\mathop{\rm sgn}(
        \eta_j)(1+\alpha - \alpha|\eta_j|)\,|\,\eta_j\big)  
        = \frac{\alpha |\eta_j|}{\alpha+1},
    \end{align}
    one easily checks that $\eta_j + \bar\zeta_j$ has the same
    distribution as $(1+\alpha)\eta_j$. Therefore, 
    $\xi_i + \zeta_i$, for $\zeta_i = a(\bar\zeta_1+\ldots + 
    \bar\zeta_k)$, has the same distribution as $(1+\alpha)
    \xi_i$. Furthermore, conditionally to the sigma-algebra 
    generated by $\{\eta_1,\ldots,\eta_k\}$, $\zeta_i$
    has zero mean and satisfies the Bernstein condition with
    $b(\alpha)=a(1+\alpha)/3$ and $v(\alpha) = a^2k\alpha 
    (1+\alpha)$.

    \begin{corollary}
        Let $\pi_0$ be a probability distribution supported by
        $S_0\subset \mathbb R^n$ such that $\mathcal D_0 = 
        \sup_{\btheta,\btheta'\in S_0}\|\btheta-\btheta'
        \|_\infty <\infty$. Assume that $\bxi$ has 
        independent entries $\xi_i$, $i=1,\ldots,n$, drawn 
        from the scaled and centered binomial distribution
        $a(\mathcal B(k,\rho_i) - k\rho_i))$. Then, for every 
        $\beta \geqslant (2/3)a\mathcal D_0$, we have    
        \begin{align}
            \bfE[\|\hat\btheta-\btheta^*\|^2]
            \leqslant \inf_{\pi}\bigg\{\int_{\mathbb R^n} 
            \|\btheta-\btheta^*\|^2 \,\pi(d\btheta) + \beta 
            D_{\rm KL}(\pi||\pi_0) \bigg\} +\bigg(\frac{6a^2k
            }{3\beta - 2a\mathcal D_0} - 1\bigg) \bfE[{\rm
            Var}_{\hat\pi}(\btheta)]. \label{eq:binomial}
        \end{align}
        In particular, if $\beta\geqslant 2a^2k + (2/3)a
        \mathcal D_0$, the last term in \eqref{eq:binomial} 
        is nonpositive and, therefore, can be neglected.
    \end{corollary}

    A typical application of this result concerns the case 
    of observing the average of $k$ Bernoulli variables, 
    that is $Y_i\sim (1/k)\mathcal B(k,\theta_i^*)$. In 
    this case, all the $\theta_i^*$ belong to $[0,1]$ and, 
    therefore, it is reasonable to choose a prior 
    distribution $\pi_0$ supported by $[0,1]^n$. This 
    ensures that $\mathcal D_0 \leqslant 1$, and, therefore, 
    inequality \eqref{eq:main1} follows from the last 
    corollary provided that $\beta\geqslant 8/(3k)$ (this 
    is obtained by choosing $a = 1/k$). 

    \subsection{Double exponential noise}

    All the previous examples considered in this section
    are distributions with sub-exponential tails. Let
    us check that \Cref{thm} can also be applied to
    some distributions that have heavier, say 
    sub-exponential, tails. Let $\xi_i$  be independent
    drawn from the Laplace distribution\footnote{This 
    means that the density of $\xi_i$ is equal to 
    $(2\mu_i)^{-1}\exp(-|x|/\mu_i)$.} with parameters $\mu_i 
    >0$, $i=1,\ldots,n$. Then, one can choose $\mathcal 
    F = \mu(\bxi)$ and $\zeta_1,\ldots, \zeta_n$ to be 
    independent, independent of $\bxi$, and drawn from 
    the distribution $\frac1{(1+\alpha )^{2}}\delta_0 
    + \frac{2\alpha + \alpha^2}{(1 + \alpha)^2} 
    \textsf{Lap}((1+\alpha)\mu_i)$. The fact that 
    $\xi_i+\zeta_i$ has the same distribution as $(1 + 
    \alpha)\xi_i$ can be checked by computing the
    characteristic functions of these variables and by 
    verifying that they are equal. As for the Bernstein
    condition, for every $t$ such that $(1+\alpha)\mu_i
    |t|\leqslant 1$ we have
    \begin{align}
        \bfE[e^{t\zeta_i}] &= \frac1{(1+\alpha)^2} + 
        \frac{2\alpha + \alpha^2}{(1+\alpha)^2} \times
        \frac{1}{1-(1+\alpha)^2 t^2\mu_i^2}\qquad 
        \big(p:=1-(1+\alpha)^{-2}, z:= (1+\alpha)t
        \mu_i\big)\\
        & = 1-p + \frac{p}{1-z^2} = 1+ \frac{pz^2}{
        1-z^2} \leqslant 1+ \frac{pz^2}{1-|z|} \\ 
        &\leqslant \exp\Big\{\frac{pz^2}{1-|z|}\Big\}
         = \exp\Big\{\frac{\alpha(2+\alpha)\mu_i^2 
         t^2}{1-(1+\alpha) \mu_i |t|}\Big\}
    \end{align}
    This means that the (conditional) Bernstein condition
    is satisfied with $v(\alpha) = \alpha(2+\alpha)
    \mu^2$ and $b(\alpha) = (1+\alpha)\mu$, where
    $\mu$ is the largest value among $\mu_i$. 

    \begin{corollary}
        Let $\pi_0$ be a probability distribution 
        supported by $S_0\subset \mathbb R^n$ such that 
        $\mathcal D_0 = \sup_{\btheta,\btheta'\in S_0} 
        \|\btheta-\btheta'\|_\infty <\infty$. Assume 
        that $\bxi$ has independent entries $\xi_i$, 
        $i=1,\ldots,n$, drawn from the Laplace 
        distribution $\textsf{Lap}(\mu_i)$. Set $
        \mu = \max_{1\leqslant i \leqslant n}\mu_i$.
        Then, for every $\beta \geqslant 2\mu
        \mathcal D_0$, we have    
        \begin{align}
            \bfE[\|\hat\btheta-\btheta^*\|^2]
            \leqslant \inf_{\pi}\bigg\{\int_{\mathbb 
            R^n} \|\btheta-\btheta^*\|^2 \,\pi 
            (d\btheta) + \beta D_{\rm KL}(\pi||\pi_0) 
            \bigg\} +\Big(\frac{4\mu^2}{\beta - 
            2\mu\mathcal D_0} - 1\Big) \bfE [ 
            {\rm Var}_{\hat\pi}(\btheta)].
            \label{eq:Laplace}
        \end{align}
        In particular, if $\beta\geqslant 4\mu^2 + 
        2\mu \mathcal D_0$, the last term in 
        \eqref{eq:Laplace} is nonpositive and, 
        therefore, can be neglected.
    \end{corollary}

    The last claim improves on 
    \citep[Prop.\ 1]{DT08}, since the latter
    requires the condition $\beta\geqslant 
    (16\mu^2) \vee (\sqrt{8}\,\mu \mathcal D_0)$. 

    \begin{remark}
        Let us finally remark that the construction
        of $\zeta_i$'s used in this section can be
        extended to the case where $\xi_i$'s are
        scale-mixtures of Laplace distributions with
        a mixing density supported by a compact set. 
        The only modification in the statement of the
        final result should be the definition of $\mu$,
        which should correspond to the smallest real
        number such that the mixing density has no mass
        in $(\mu,\infty)$. Similar extension can be 
        carried out in the case of scale-mixtures 
        of Gaussians.
    \end{remark}


\section{Proof of \Cref{thm}}\label{sec:proof}

    Since $\hat\pi$ minimizes the criterion $\pi\mapsto 
    \int_{\mathbb R^n} \|\bY-\btheta\|^2\,\pi(d\btheta) + 
    \beta D_{\rm KL}(\pi||\pi_0)$, we have
    \begin{align}
        \int_{\mathbb R^n} \|\bY-\btheta\|^2\,\hat\pi
        (d\btheta) + \beta D_{\rm KL}(\hat\pi||\pi_0) 
        \leqslant \int_{\mathbb R^n} \|\bY-\btheta\|^2 
        \,\pi(d\btheta) + \beta D_{\rm KL}(\pi||\pi_0)
        \label{eq:2}
    \end{align}
    for all densities $\pi$ over $\mathbb R^n$. The 
    KL-divergence being always nonnegative, we infer from 
    the last display that
    \begin{align}
        \|\bY-\hat\btheta\|^2 & = \int_{\mathbb R^n} \|\bY - 
        \btheta\|^2\,\hat\pi(d\btheta) - \int_{\mathbb R^n}
        \|\btheta - \hat\btheta\|^2\,\hat\pi(d\btheta)\\
        &\leqslant \int_{\mathbb R^n} \|\bY-\btheta\|^2 
        \,\pi(d\btheta) + \beta D_{\rm KL}(\pi||\pi_0)
        - \int_{\mathbb R^n} \|\btheta  -  \hat\btheta\|^2 
        \, \hat\pi(d\btheta).\label{eq:3}
    \end{align}
    Using the decompositions $\|\bY-\hat\btheta\|^2 = 
    \|\hat\btheta-\btheta^*\|^2 +2(\hat\btheta-\btheta^*)^\top\bxi
    +\|\bxi\|^2$ and $\|\bY-\btheta\|^2 = \|\btheta-\btheta^*\|^2
    +2(\btheta^* - \btheta)^\top\bxi + \|\bxi\|^2$ and taking the 
    expectation of the two sides of \eqref{eq:3}, we get
    \begin{align}
        \bfE[\|\hat\btheta-\btheta^*\|^2] + 
        2\bfE[(\hat\btheta-\btheta^*)^\top\bxi] 
        \leqslant \int_{\mathbb R^n} \|\btheta-\btheta^*\|^2 
        \,\pi(d\btheta) + \beta D_{\rm KL}(\pi||\pi_0) - 
        \bfE\bigg[\int_{\mathbb R^n} \|\btheta  - \hat\btheta\|^2 
        \, \hat\pi(d\btheta)\bigg]
    \end{align}
    which can be equivalently written as
    \begin{align}
        \bfE[\|\hat\btheta-\btheta^*\|^2]
        \leqslant \int_{\mathbb R^n} \|\btheta-\btheta^*\|^2 
        \,\pi(d\btheta) + \beta D_{\rm KL}(\pi||\pi_0)
        + 2\bfE[\hat\btheta{}^\top\bxi] 
        - \int_{\mathbb R^n} \bfE[\|\btheta  -  \hat\btheta\|^2 
        \, \hat\pi(\btheta)]\,d\btheta.\label{eq:1}
    \end{align}
    In addition, we have
    \begin{align}
        2\bfE[\hat\btheta{}^\top\bxi]  = \frac\beta{\alpha}
        \bfE\bigg[\int_{\mathbb R^n} \log e^{2(\alpha/\beta)
        \btheta^\top \bxi}\hat\pi(d\btheta)\bigg],
    \end{align}
    where $\alpha>0$ is an arbitrary number. Since the logarithm
    is concave, the Jensen inequality yields
    \begin{align}
        2\bfE[\hat\btheta{}^\top\bxi]  &\leqslant \frac\beta{
        \alpha} \bfE \bigg[ \log\bigg(\int_{\mathbb R^n}
        e^{2(\alpha/\beta) \btheta^\top \bxi} \hat\pi(d\btheta)
        \bigg) \bigg] \\ 
        &=   \frac\beta{\alpha} \bfE \bigg[\log\bigg(\int_{
        \mathbb R^n} e^{2(\alpha/\beta) \btheta^\top \bxi - 
        \|\btheta^*+\bxi - \btheta\|^2/\beta} \, \pi_0 
        (d\btheta)\bigg)  - \log\bigg(\int_{\mathbb R^n} 
        e^{- \|\btheta^*+\bxi - \btheta\|^2/\beta} \, 
        \pi_0(d\btheta) \bigg) \bigg]\\
        &=   \frac\beta{\alpha} \bfE \bigg[ \log\bigg(\int_{
        \mathbb R^n} e^{(2(1+\alpha)\btheta^\top \bxi 
        -  \|\btheta^* - \btheta\|^2)/\beta} \,\pi_0 
        (d\btheta)\bigg) - \log\bigg(\int_{\mathbb R^n} 
        e^{(2 \btheta^\top \bxi- \|\btheta^* - \btheta\|^2)
        /\beta} \, \pi_0(d\btheta)\bigg) \bigg]\label{eq:4}
    \end{align}
    Let $\bzeta = \bzeta_\alpha$ be the $n$ dimensional random
    vector the existence of which is required in the statement
    of the theorem. Recall that it satisfies 
    \begin{align}\label{cond:zeta}
        \bfE[\bzeta|\mathcal F] = 0,\qquad \bxi + \bzeta 
        \stackrel{\mathscr D}{=} (1+\alpha) \bxi ,\qquad 
    \end{align}
    These conditions imply that in the first expectation
    in \eqref{eq:4}, one can replace $(1+\alpha)\bxi$ 
    by $\bxi+\bzeta$, which yields
    \begin{align}
        2\bfE[\hat\btheta{}^\top\bxi]  &\leqslant 
        \frac\beta{\alpha} \bfE \bigg[ \log\bigg(\int_{
        \mathbb R^n} e^{(2\btheta^\top\bxi + 2\btheta^\top
        \bzeta - \|\btheta^* - \btheta\|^2)/\beta} \,
        \pi_0(d\btheta)\bigg) \bigg] - \frac\beta{\alpha} 
        \bfE \bigg[ \log\bigg(\int_{\mathbb R^n} e^{(2
        \btheta^\top \bxi- \|\btheta^* - \btheta\|^2)/\beta} \,\pi_0(d\btheta)\bigg) \bigg]\\
        & = \frac\beta{\alpha} \bfE \bigg[ \log\bigg( \int_{
        \mathbb R^n} e^{2\btheta^\top\bzeta/\beta} \,\hat\pi 
        (d\btheta)\bigg) \bigg]
         = \frac\beta{\alpha} \bfE \bigg[ \log\bigg( \int_{
        \mathbb R^n} e^{2(\btheta-\hat\btheta)^\top\bzeta/\beta}
        \,\hat\pi(d\btheta)\bigg) \bigg].
        \label{eq:5}
    \end{align}
    Since conditionally to $\bxi$, $\zeta_i$'s are independent
    and each $\zeta_i$ satisfies the $(v(\alpha), b(\alpha)) 
    $-Bernstein condition, one can use the Jensen inequality 
    to upper bound the expectation in\eqref{eq:5} as follows
    \begin{align}
        2\bfE[\hat\btheta{}^\top\bxi]  & \leqslant
        \frac\beta{\alpha} \bfE \bigg[ \log\bigg( \int_{
        \mathbb R^n} \bfE[e^{2(\btheta-\hat\btheta)^\top\bzeta
        /\beta} |\mathcal F]\,\hat\pi (d\btheta)\bigg) \bigg]\\
        &\leqslant \frac\beta{\alpha} \bfE \bigg[ 
        \log\bigg( \int_{ \mathbb R^n} \exp\Big\{\frac{2
        \|\btheta - \hat\btheta\|^2 v(\alpha)}{\beta(\beta 
        - 2 b(\alpha) \|\btheta-\hat\btheta\|_\infty)}\Big\}
        \,\hat\pi (d\btheta) \bigg) \bigg]\label{eq:6}
    \end{align}
    for every $\beta$ satisfying $\beta\geqslant 2b(\alpha)
    \|\btheta - \btheta'\|_\infty$ for every $\btheta,\btheta'
    \in S_0 : = \text{supp}(\pi_0)$. Note that for every 
    $\btheta\in S_0$, we have $\|\btheta-\hat\btheta\|_\infty
    \leqslant \mathcal D_\infty(S_0)$. The inequality in 
    \eqref{eq:6} being true for any $\alpha>0$, one can 
    check that 
    \begin{align}
        2\bfE[\hat\btheta{}^\top\bxi]  
        &\leqslant \liminf_{\alpha\to 0}\frac\beta{\alpha} 
        \bfE \bigg[ \log\bigg( \int_{ \mathbb R^n} 
        \exp\Big\{\frac{2\|\btheta - \hat\btheta\|^2 
        v(\alpha)}{\beta( \beta - 2 b(\alpha)\|\btheta -
        \hat\btheta\|_\infty)}\Big\}\,\hat\pi (d\btheta)
        \bigg) \bigg]\\
        & = \bfE \bigg[ \int_{ \mathbb R^n} \frac{2
        \|\btheta - \hat\btheta\|^2 v'(0)}{\beta - 2b(0)
        \|\btheta - \hat\btheta\|_\infty}\,\hat\pi 
        (d\btheta) \bigg]
        \leqslant \frac{2 v'(0)}{\beta - 2 b(0) 
        \mathcal D_\infty(S_0)}\bfE \bigg[\int_{\mathbb R^n} 
        \|\btheta - \hat\btheta\|^2\,\hat\pi (d\btheta) 
        \bigg]. \label{eq:7}
    \end{align}
    Combining \eqref{eq:1} and \eqref{eq:7}, we see that 
    \begin{align}
        \bfE[\|\hat\btheta-\btheta^*\|^2]
        \leqslant \int_{\mathbb R^n} \|\btheta-\btheta^*\|^2 
        \,\pi(d\btheta) + \beta D_{\rm KL}(\pi||\pi_0)
        +\bigg(\frac{2 v'(0)}{\beta - 2 b(0) \mathcal
        D_\infty(S_0)}- 1\bigg) \bfE[\text{Var}_{\hat\pi}
        (\btheta)].\label{eq:8}
    \end{align}
    This completes the proof. 

    \paragraph*{Acknowledgments}  The work of the author 
    was supported by the grant Investissements d’Avenir 
    (ANR-11-IDEX-0003/Labex Ecodec/ANR-11-LABX-0047),  
    the FAST Advance grant and the center Hi! PARIS. 
   \bibliography{bibliography}

\begin{thebibliography}{22}
\providecommand{\natexlab}[1]{#1}
\providecommand{\url}[1]{\texttt{#1}}
\expandafter\ifx\csname urlstyle\endcsname\relax
  \providecommand{\doi}[1]{doi: #1}\else
  \providecommand{\doi}{doi: \begingroup \urlstyle{rm}\Url}\fi

\bibitem[Alquier and Lounici(2011)]{Alquier3}
Pierre Alquier and Karim Lounici.
\newblock P{AC}-{B}ayesian bounds for sparse regression estimation with
  exponential weights.
\newblock \emph{Electron. J. Stat.}, 5:\penalty0 127--145, 2011.

\bibitem[Bellec(2018)]{bellec2018}
Pierre~C. Bellec.
\newblock Optimal bounds for aggregation of affine estimators.
\newblock \emph{Ann. Statist.}, 46\penalty0 (1):\penalty0 30--59, 02 2018.
\newblock \doi{10.1214/17-AOS1540}.

\bibitem[Bunea et~al.(2007)Bunea, Tsybakov, and Wegkamp]{BTW2007}
Florentina Bunea, Alexandre~B. Tsybakov, and Marten~H. Wegkamp.
\newblock Aggregation for gaussian regression.
\newblock \emph{Ann. Statist.}, 35\penalty0 (4):\penalty0 1674--1697, 08 2007.

\bibitem[Chernousova et~al.(2013)Chernousova, Golubev, and Krymova]{Golubev1}
Elena Chernousova, Yuri Golubev, and Ekaterina Krymova.
\newblock Ordered smoothers with exponential weighting.
\newblock \emph{Electron. J. Stat.}, 7:\penalty0 2395--2419, 2013.

\bibitem[Dalalyan and Tsybakov(2009)]{colt_DalalyanT09}
A.~S. Dalalyan and A.~B. Tsybakov.
\newblock Sparse regression learning by aggregation and {L}angevin
  {M}onte-{C}arlo.
\newblock In \emph{COLT 2009 - The 22nd Conference on Learning Theory,
  Montreal, June 18-21, 2009}, pages 1--10, 2009.

\bibitem[Dalalyan(2020)]{Dal_IHP}
Arnak~S. Dalalyan.
\newblock {Exponential weights in multivariate regression and a low-rankness
  favoring prior}.
\newblock \emph{Annales de l'Institut Henri Poincaré, Probabilités et
  Statistiques}, 56\penalty0 (2):\penalty0 1465 -- 1483, 2020.
\newblock \doi{10.1214/19-AIHP1010}.
\newblock URL \url{https://doi.org/10.1214/19-AIHP1010}.

\bibitem[Dalalyan and Salmon(2012)]{DalSal}
Arnak~S. Dalalyan and Joseph Salmon.
\newblock {Sharp oracle inequalities for aggregation of affine estimators}.
\newblock \emph{The Annals of Statistics}, 40\penalty0 (4):\penalty0 2327 --
  2355, 2012.
\newblock \doi{10.1214/12-AOS1038}.
\newblock URL \url{https://doi.org/10.1214/12-AOS1038}.

\bibitem[Dalalyan and Tsybakov(2007)]{DT07}
Arnak~S. Dalalyan and Alexandre~B. Tsybakov.
\newblock Aggregation by exponential weighting and sharp oracle inequalities.
\newblock In \emph{Learning theory}, volume 4539 of \emph{Lecture Notes in
  Comput. Sci.}, pages 97--111. Springer, Berlin, 2007.

\bibitem[Dalalyan and Tsybakov(2008)]{DT08}
Arnak~S. Dalalyan and Alexandre~B. Tsybakov.
\newblock Aggregation by exponential weighting, sharp pac-bayesian bounds and
  sparsity.
\newblock \emph{Machine Learning}, 72\penalty0 (1-2):\penalty0 39--61, 2008.

\bibitem[Donier-Meroz et~al.(2023)Donier-Meroz, Dalalyan, Kramarz, Chon{\'e},
  and {D'Haultfoeuille}]{Etienne1}
E.~Donier-Meroz, A.~S. Dalalyan, F.~Kramarz, Ph. Chon{\'e}, and
  X.~{D'Haultfoeuille}.
\newblock Graphon estimation in bipartite graphs with observable edge labels
  and unobservable node labels.
\newblock Technical report, 2023.

\bibitem[George(1986{\natexlab{a}})]{George86b}
E.~I. George.
\newblock Combining minimax shrinkage estimators.
\newblock \emph{J. Amer. Statist. Assoc.}, 81\penalty0 (394):\penalty0
  437--445, 1986{\natexlab{a}}.

\bibitem[George(1986{\natexlab{b}})]{George86}
Edward~I. George.
\newblock {Minimax Multiple Shrinkage Estimation}.
\newblock \emph{The Annals of Statistics}, 14\penalty0 (1):\penalty0 188 --
  205, 1986{\natexlab{b}}.
\newblock \doi{10.1214/aos/1176349849}.
\newblock URL \url{https://doi.org/10.1214/aos/1176349849}.

\bibitem[Juditsky et~al.(2008)Juditsky, Rigollet, and Tsybakov]{JRT}
A.~Juditsky, P.~Rigollet, and A.~B. Tsybakov.
\newblock Learning by mirror averaging.
\newblock \emph{Ann. Statist.}, 36\penalty0 (5):\penalty0 2183--2206, 2008.

\bibitem[Lecu{\'e} and Mendelson(2013)]{Lecue}
Guillaume Lecu{\'e} and Shahar Mendelson.
\newblock On the optimality of the aggregate with exponential weights for low
  temperatures.
\newblock \emph{Bernoulli}, 19\penalty0 (2):\penalty0 646--675, 2013.

\bibitem[Leung and Barron(2006)]{LeungBarron}
G.~Leung and A.R. Barron.
\newblock Information theory and mixing least-squares regressions.
\newblock \emph{IEEE Transactions on Information Theory}, 52\penalty0
  (8):\penalty0 3396--3410, 2006.
\newblock \doi{10.1109/TIT.2006.878172}.

\bibitem[Rigollet and Tsybakov(2011)]{Rigollet11}
Philippe Rigollet and Alexandre Tsybakov.
\newblock Exponential screening and optimal rates of sparse estimation.
\newblock \emph{Ann. Statist.}, 39\penalty0 (2):\penalty0 731--771, 2011.

\bibitem[Tsybakov(2003)]{Tsyb03}
Alexandre~B. Tsybakov.
\newblock Optimal rates of aggregation.
\newblock In Bernhard Sch{\"o}lkopf and Manfred~K. Warmuth, editors,
  \emph{Learning Theory and Kernel Machines}, pages 303--313, Berlin,
  Heidelberg, 2003. Springer Berlin Heidelberg.

\bibitem[Tsybakov(2014)]{TsybICM}
Alexandre~B. Tsybakov.
\newblock Aggregation and minimax optimality in high-dimensional estimation.
\newblock In \emph{Proceedings of the International Congress of Mathematicians
  (Seoul, August 2014)}, volume~3, pages 225--246, 2014.

\bibitem[Vershynin(2018)]{vershynin_2018}
Roman Vershynin.
\newblock \emph{High-Dimensional Probability: An Introduction with Applications
  in Data Science}.
\newblock Cambridge Series in Statistical and Probabilistic Mathematics.
  Cambridge University Press, 2018.
\newblock \doi{10.1017/9781108231596}.

\bibitem[Yang(2000)]{Yang00}
Y.~Yang.
\newblock Combining different procedures for adaptive regression.
\newblock \emph{J. Multivariate Anal.}, 74\penalty0 (1):\penalty0 135--161,
  2000.

\bibitem[Yang(2003)]{Yang2}
Yuhong Yang.
\newblock Regression with multiple candidate models: selecting or mixing?
\newblock \emph{Statist. Sinica}, 13\penalty0 (3):\penalty0 783--809, 2003.

\bibitem[Yang(2004)]{Yang1}
Yuhong Yang.
\newblock Aggregating regression procedures to improve performance.
\newblock \emph{Bernoulli}, 10\penalty0 (1):\penalty0 25--47, 2004.

\end{thebibliography}
\end{document}